\documentclass[12pt]{amsart}
\newtheorem{theorem}{Theorem}

\newtheorem{lemma}{Lemma}
\newtheorem{proposition}{Proposition}
\newtheorem{definition}{Definition}
\newtheorem{corollary}{Corollary}
\newtheorem{example}{Example}

\newenvironment{pf}{\noindent {\bf Proof:}}{{\qed}}

\newcommand{\df}{:=}
\newcommand{\Nnn}{{\mathbb N}}
\newcommand{\Ppp}{{\mathbb P}}

\newcommand{\Rrr}{{\mathbb R}}

\newcommand{\vanish}[1]{}
\newcommand{\iv}{{\mathcal I}}
\newcommand{\rank}{\mbox{\rm rank}}
\newcommand{\1}{\hat{1}}
\newcommand{\0}{\hat{0}}
\newcommand{\block}{\mbox{\rm\bf B}}
\newcommand{\cone}[1]{K_{#1}}
\newcommand{\flag}[1]{A_{#1}}
\newcommand{\dflag}[1]{A_{#1}^{\ast}}
\newcommand{\proj}{{\bf \pi}}
\newcommand{\cut}{{\bf \rho}}
\newcommand{\supp}{\mbox{\rm supp}}
\newcommand{\faa}{K\langle x_1,x_2, \ldots\rangle}
\newcommand{\eval}[1]{\varepsilon_{#1}}
\newcommand{\KER}{\mbox{\rm Ker}}
\newcommand{\IM}{\mbox{\rm Im}}
\newcommand{\shift}[1]{\sigma_{#1}}

\begin{document}

\title[Flags in Graded Posets]{Linear inequalities for flags in graded 
posets$^{\ast}$}\thanks{$^{\ast}$Research at MSRI supported in part by NSF grant
DMS-9022140.}

\author[L.\ J.\ Billera]{Louis J.\ Billera$^{\dagger}$}\thanks{$^{\dagger}$On leave from 
	Cornell University, Ithaca, NY.  Supported in part by
	NSF grant DMS-9500581.}
\author[G.\ Hetyei]{G\'abor Hetyei$^{\ddagger}$}\thanks{$^{\ddagger}$On leave from the Mathematical Research
        Institute of the Hungarian Academy of Sciences.
	Partially supported by Hungarian National Foundation for
	Scientific Research grant no. F 023436}
\address{	MSRI\\
	1000 Centennial Drive\\
	Berkeley, CA 94720-5070}

\begin{abstract}
The closure of the convex cone generated by all flag $f$-vectors of
graded posets is shown to be polyhedral.  In particular, we give
the facet inequalities to the polar
cone of all nonnegative chain-enumeration functionals on this class of
posets.  These are in one-to-one correspondence with antichains of intervals
on the set of ranks and thus are counted by Catalan
numbers. Furthermore, we prove that the convolution operation introduced
by Kalai assigns extreme rays to pairs of extreme rays in most cases. We
describe the strongest possible inequalities for graded posets of rank
at most 5.
\end{abstract}

\maketitle

\section*{Introduction}
An initial step in obtaining a characterization of $f$-vectors of some
class of objects is to determine the linear equations and inequalities that
they must satisfy.  The former give a description of the linear span of
all such $f$-vectors, while the latter describe the closure of the convex
cone they generate.  In most cases where this has been done successfully,
the description of the linear equations proved to be the more difficult part.
Once this was done, and an appropriate basis found for the linear span of
all $f$-vectors (for example, the $f$-vector in the case of simplicial
complexes, the $h$-vector in the case of
Cohen-Macaulay simplicial complexes, or the $g$-vector in the case of
simplicial convex polytopes), the desired cone turned out to be
an orthant, that is, only nonnegativity of the basic invariants could
be asserted.  (See, for example, \cite[Theorems II.2.1, II.3.3 \& 
III.1.1]{CCA}.)

The situation of flag $f$-vectors seems to be be quite different.
While it is true that flag $h$-vectors of balanced Cohen-Macaulay
complexes (even Cohen-Macaulay ranked posets -- see \cite[\S III.4]{CCA})
span an orthant, the more basic case of flag $f$-vectors of graded posets
behaves quite differently.  Here, it is
the first of these that is simple: there are no equations.  However, the
cone generated in this case turns out to be quite a bit more
complicated.  In this paper, we
give a description of this cone by giving its minimal
generating set.  It is already a nontrivial statement that this
generating set is finite.  Equivalently, we give a finite list of linear
inequalities that describe the cone polar to that generated by
all flag $f$-vectors
of graded posets.  For posets of rank $n+1$, these inequalities are in
one-one correspondence with antichains of intervals in the linearly
ordered set $\{1,2,\dots,n\}$, an so are counted by a Catalan number.
Thus while the space of flag $f$-vectors has dimension $2^n$,
the cone they generate will have on the order of $(2n)^n$ generators.

The proof of the fact that our list of linear inequalities describes a
set of generators for the cone of flag $f$-vectors relies on two
ingredients: an explicit construction of sequences of graded posets
$P(n,\iv,N)$ of rank $n+1$ yielding the extreme rays of the closure as
limits, and an explicit partitioning of the maximal chains of every
graded poset $P$ which allows to show the sufficiency of our list of
conditions. The chain partitioning used may be generalized to a
construction showing that every graded poset satisfies a generalized
condition of lexicographic shellability, allowing an explicit
description of the order complex of every graded poset. This will be the
subject of a subsequent paper. 

While the description of the extreme rays of the closure of the convex
cone generated by the flag $f$-vectors of graded posets of a given rank
is fairly tractable, finding even the number of the facets of the same
cone seems to be highly difficult. These represent the strongest
possible inequalities holding for flag $f$-vectors of graded posets of a
given rank, and correspond to the extreme rays of the polar cone. In
section \ref{S_extreme} we show operations which yield higher-rank
extremes from lower rank ones. There are lifting operations, embedding
the cones of inequalities into each other as faces, and we give an exact
description of those situation where the convolution operation
introduced by Kalai in \cite{Kalai} assigns an extreme inequality to a
pair of extreme inequalities. These results allow a short description of
the extreme inequalities up to rank $5$. 

\section{Preliminaries}
Here we enumerate the basic definitions and results used in this paper.
\subsection{Graded partially ordered sets}
\label{ss_gpo}
\begin{definition}
A {\em graded} poset $P$ is a finite poset with a unique minimum element
$\0$, a unique maximum element $\1$, and a {\em rank function} $\rank:
P\longrightarrow \Nnn$ such that we have 
\begin{enumerate}
\item[(i)] $\rank (\0)=0$, and
\item[(ii)] $\rank(y)-\rank(x)=1$ whenever $y\in P$ covers $x\in P$.
\end{enumerate}
\end{definition}

We call $\rank(\1) \ge 1$ the rank of the poset $P$. Given a graded poset $P$
of rank $n+1$, and a subset $S$ of $\{1,2,\ldots,n\}$ we define the {\em
$S$--rank selected subposet of $P$} to be the poset
$$P_{S} \df \{ x \in P\::\: \rank(x) \in S\} \cup \{ {\0},{\1}\}.$$
We denote by $f_S (P)$ the number of maximal chains of $P_S$. Equivalently,
$f_S(P)$ is the number of chains $x_1<\cdots<x_{|S|}$ in $P$ such that 
$\{\rank(x_1),\ldots,\rank(x_{|S|})\}=S$. The function 
$$\begin{array}{rrcl}
f:& 2^{\{1,2,\ldots,n\}}&\longrightarrow& \Nnn\\
&S&\longmapsto& f_S (P)\\
\end{array}$$ 
is called the {\em flag
$f$-vector} of $P$. Whenever it does not cause confusion we will write
$f_{s_1\,\ldots\, s_k}$ rather than $f_{\{s_1,\ldots,s_k\}}$; in particular,
$f_{\{m\}}$ will always be denoted $f_m$. 

A simple example is perhaps helpful here.  Consider graded posets of
rank 2.  In addition to elements $\hat 0$ and $\hat 1$ of rank 0 and 2,
respectively, such a poset will have only elements of rank 1.
For reasons that will become clear later, we wish to consider the closure
of the convex cone generated by
$$\{(f_\emptyset(P),f_1(P)) \ | \ \rank(P)=2 \}\ =\ \{ (1,n) \ | \ n\ge 0 \}.$$
The polar of this cone is generated by the functionals 
$h_\emptyset \df f_\emptyset$ and $h_1 \df f_1 - f_\emptyset$.

\subsection{The ring of chain operators}

We adopt the following terminology and results from
\cite{Billera-Liu}. The {\em chain operators}  $f_S^{n+1}$ assign
$f_S(P)$ to every graded poset $P$ of rank $n+1$, and zero to all other
graded posets. These operators are linearly independent and hence they
generate a vector space $\flag{n+1}$ of dimension $2^n$ over $\Rrr$. We set
$\flag{0}=\Rrr$. The vector space $\flag{}\df \bigoplus_{n\geq 0}
\flag{n}$ may be made into a graded noncommutative ring by introducing the {\em
convolution operation} (first considered by Kalai in \cite{Kalai})
$$f^m_S \ast f^n _T\df f^{m+n}_{S\cup \{m\}\cup (T+m)}$$ 
for $m,n\geq 1$, and by making the generator $1$ of $\flag{0}$ to be the
unit of $\flag{}$.   
The interest of this convolution operation is the
following. (cf. \cite[Proposition 1.3]{Billera-Liu}) 
\begin{proposition}
\label{P_conv}
The convolution $F\ast G$ of two linear combinations of chain operators
is nonnegative on all graded posets if and only if both $F$ and $G$ are 
simultaneously nonnegative or nonpositive on all graded posets.
\end{proposition}
According to \cite[Proposition 1.1]{Billera-Liu} the
chain operators are linearly independent.  Moreover, in \cite[Section
2]{Billera-Liu} we find the following: 
\begin{proposition}
\label{P_free}
The ring $\flag{}$ is a free graded associative algebra over the set of
variables  $\left\{f^n_{\emptyset}\::\: n\geq 1 \right\}$.
\end{proposition}
It follows from \cite[Theorem 3]{Cohn}
that the semigroup of homogeneous
polynomials of a a free graded associative
algebra has unique factorization. Hence we have the following
\begin{proposition}
\label{P_ufact}
Up to nonzero linear factors, every form $F\in \flag{n}$ may be uniquely
written as a product $F=F^{n_1}\ast\cdots\ast F^{n_k}$ of homogeneous forms
$F^{n_i}\in\flag{n_i}$ such that we have $n=n_1+\cdots+n_k$ and the
$F^{n_i}$'s can not be written as a product of two homogeneous forms of
lower degree.     
\end{proposition}

\subsection{Blockers of families of sets}

\begin{definition}
Let ${\mathcal S}$ be an arbitrary family of subsets of a finite set $X$. A subset
$B\subseteq X$ is a {\em blocker} of ${\mathcal S}$, if for
every $S\in {\mathcal S}$ we have $S\cap T\neq\emptyset$.
\end{definition}
We denote the set of blockers of ${\mathcal S}$ in $X$  by $\block _X ({\mathcal
S})$. In particular, for ${\mathcal S}=\emptyset$ every subset of $X$ is a
blocker and so we have $\block _X({\emptyset})=\{T\::\:T\subseteq X\}$.   

\begin{lemma}
Let ${\mathcal S}_1$ and ${\mathcal S}_2$ be a family of subsets of the same
finite set $X$. Then $\block _X ({\mathcal S}_1)\subseteq\block _X ({\mathcal
S}_2)$ if and only if every $S_2\in {\mathcal S}_2$ contains some $S_1\in
{\mathcal S}_1$.  
\label{L_block}
\end{lemma}
\begin{proof}
If every set from ${\mathcal S}_2$ contains a set from ${\mathcal S}_1$ then
every blocker of ${\mathcal S}_1$ also blocks ${\mathcal S}_2$. Assume that
there is an element $S_2$ of ${\mathcal S}_2$ not containing any element of ${\mathcal
S}_1$. Then the set $X\setminus S_2$ blocks ${\mathcal S}_1$, but it is
disjoint from $S_2$ and so does not block ${\mathcal S}_2$. Hence we have
$\block _X ({\mathcal S}_1)\not \subseteq \block _X ({\mathcal S}_2)$. 
\end{proof}

\begin{corollary}
\label{C_block}
Let ${\mathcal S}_1$ and ${\mathcal S}_2$ be a family of subsets of the same
finite set $X$. Then $\block _X ({\mathcal S}_1)=\block _X ({\mathcal S}_2)$ 
if and only if the minimal sets with respect to inclusion are the same
in ${\mathcal S}_1$ and ${\mathcal S}_2$. 
\end{corollary}

The minimal sets in a family of sets form an {\em antichain} (or
Sperner family, or clutter), i.e. a family
of sets such that no set contains another one. On the other hand, a
family of sets $\mathcal S$ on a set $X$  
is a {\em dual ideal} if for every $S\in {\mathcal S}$ all $S'\subseteq X$
containing $S$ belong to $\mathcal S$.
For an arbitrary family of sets $\mathcal S$ on $X$ the {\em dual
ideal generated by $\mathcal S$} is the family
$${\mathcal S} ^+\df \{T\subseteq X\::\: \exists S\in {\mathcal S}
(S\subseteq T)\}.$$ 

\begin{proposition}
\label{P_dual}
For any family of sets $\mathcal S$ on a set $X$, we
have $\block _X({\mathcal S})=\block _X({\mathcal S}^+)$.
\end{proposition}
\begin{proof}
Since $\mathcal S$ is a subfamily of ${\mathcal S}^+$, we have $\block
_X({\mathcal S})\supseteq\block _X({\mathcal S} ^+)$. The other inclusion follows from
Lemma \ref{L_block} using the fact that every set from  ${\mathcal S} ^+$
contains some set from $\mathcal S$.
\end{proof}

It is shown in
\cite{Edmonds-Fulkerson} that whenever ${\mathcal S}$ is an antichain,
then the family of minimal elements of $\block _X(\block _X({\mathcal
S}))$ is ${\mathcal S}$. This statement and repeated
application of Proposition \ref{P_dual} yields
$$\block _X(\block _X({\mathcal S}))={\mathcal S}^+$$
for every family of sets $\mathcal S$. From this and Proposition
\ref{P_dual} we have the following. 

\begin{proposition}
\label{P_dual2}
Let ${\mathcal S}_1$ and ${\mathcal S}_2$ be families of sets on the same set
$X$. Then 
$$\block _X({\mathcal S}_1)=\block _X({\mathcal S}_2)\mbox{ if and
only if } {\mathcal S}_1^+={\mathcal S}_2^+.$$
\end{proposition} 

For more on blockers, we refer the reader to
\cite{Edmonds-Fulkerson}, or to Section 8.1 of
\cite{Grotschel-Lovasz-Schrijver}. 

\section{The Main Theorem}
\label{S_main}

This section contains our main result, which characterizes all linear 
inequalities holding for the flag $f$-vectors of all graded posets of
rank $n+1$. Every such linear inequality may be written as 
$$F(P)\df \sum_{S\subseteq\{1,\ldots,n\}} a_S\cdot f^{n+1}_S (P)\geq 0,$$
where the coefficients $a_S$ are real numbers. Moreover,
$\sum_{S\subseteq\{1,\ldots,n\}} a_S\cdot f^{n+1}_S (P)$ is zero if
$\rank(P)\neq n+1$. Hence we are interested
in determining the subset 
$$\cone{n+1}\df \left\{ F=\sum_{S\subseteq\{1,\ldots,n\}} a_S\cdot
f^{n+1}_S\in\flag{n+1}\::\: \forall P \left(F (P)\geq 0\right)\right\}$$ 
of $\flag{n+1}$.  For convenience, we will let $\cone{0} = \flag{0}$.

In this section we show that $\cone{n+1}$ is a polyhedral cone, that is, the
intersection of finitely many half spaces. We give these half spaces in
terms of {\em interval systems} on the linearly ordered set
$\{1,2,\ldots,n\}~\subset~\Nnn$.

A subset of a partially ordered set $P$ is an {\em interval}, if it is
empty, or of the form 
$$[p,q]\df \{x\in P\::\: p\leq x\leq q\}\quad \mbox{for some $p,q\in P$.}$$ 
In particular, an {\em interval} in $\Nnn$ is a finite (possibly empty) set of
consecutive natural numbers. An {\em interval system on a partially
ordered set $P$} is a family $\iv$ of intervals. We consider the empty
set also as an interval system. The following theorem gives the list of
inequalities that determine $\flag{n+1}$. 

\begin{theorem}
\label{T_main}
An expression $\sum _{S\subseteq [1,n]}
a_S\cdot f^{n+1}_S$ is nonnegative on all graded posets of rank $n+1$  if and only
if we have
\begin{equation}
\label{E_main}
\sum _{S\in \block_{[1,n]}(\iv)} a_S\geq 0 \quad \mbox{for every interval system $\iv$ on $[1,n]$.}
\end{equation}
\end{theorem}
First we will show the necessity of condition (\ref{E_main}) by constructing
for every interval system $\iv$ on $\{1,2,\ldots,n \}$ a family of
posets $\{P(n,\iv,N)\::\: N\in \Nnn\}$  of rank $n+1$ such that we have
$$\lim _{N\longrightarrow \infty} \frac{1}{
{f_{[1,n]}(P(n,\iv,N))}}\cdot \sum _{S\subseteq
[1,n]}a_S\cdot f_S(P(n,\iv,N)))=\sum _{{S\in \block_{[1,n]}(\iv)}} a_S.$$
Then we will prove the sufficiency by using an appropriate partitioning
of the set of maximal chains for every graded poset $P$. 

\begin{definition}
Let $n$ and $N$ be positive integers, and let
$\iv=\{I_1,I_2\ldots,I_k\}$ be an interval system on $[1,n]$. We define the elements of the poset $P(n,\iv,N)$ to be 
all arrays  $(i;p_1,\ldots,p_k)$ such that
\begin{enumerate}
\item[(i)] $i\in [0,n+1]$, and 
\item[(ii)] for every $j\in[1,k]$ we have 
$$p_j\in\left\{\begin{array}{ll}
[1,N] &\mbox{whenever $i\in I_j$,}\\
\{\ast\} & \mbox{otherwise}\end{array}\right.$$
Here $\ast$ is a special symbol, different from all integers.
\end{enumerate}
We set $(i;p_1,p_2,\ldots,p_k)\leq
(i';p'_1,p'_2,\ldots,p'_k)$ if 
\begin{enumerate}
\item[(1)] $i\leq i'$, and  
\item[(2)] for every $j\in[1,k]$ we have either $p_j=p'_j$ or
$\ast\in \{p_j,p'_j\}$. 
\end{enumerate}
\end{definition}
Observe that $P(n,\iv,N)$ has a unique minimum element
$\0=(0;\ast,\ldots,\ast)$ and a unique maximum element $\1=(n+1;\ast,\ldots,\ast)$.

\begin{example}{\em
Let $n=3$, $\iv=\{[1,2],[2,3]\}$, and $N=2$. Then $P(3,\{[1,2],[2,3]\},2)$
is the poset shown on Fig. \ref{F_pin1}.  

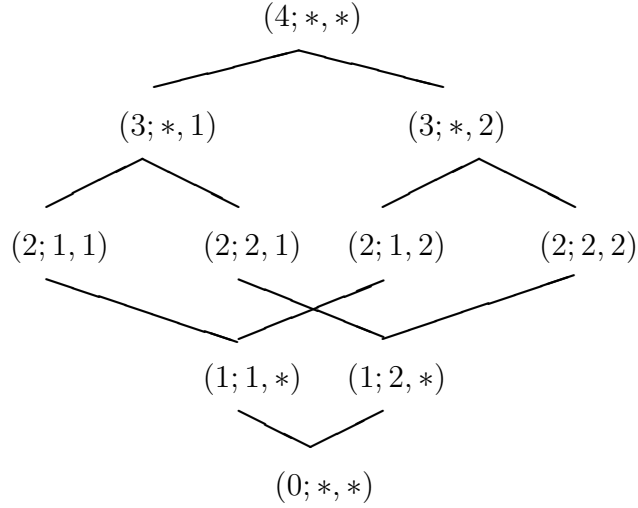
\begin{figure}
\begingroup\makeatletter
\def\x#1#2#3#4#5#6#7\relax{\def\x{#1#2#3#4#5#6}}%
\expandafter\x\fmtname xxxxxx\relax \def\y{splain}%
\ifx\x\y   
\gdef\SetFigFont#1#2#3{%
  \ifnum #1<17\tiny\else \ifnum #1<20\small\else
  \ifnum #1<24\normalsize\else \ifnum #1<29\large\else
  \ifnum #1<34\Large\else \ifnum #1<41\LARGE\else
     \huge\fi\fi\fi\fi\fi\fi
  \csname #3\endcsname}%
\else
\gdef\SetFigFont#1#2#3{\begingroup
  \count@#1\relax \ifnum 25<\count@\count@25\fi
  \def\x{\endgroup\@setsize\SetFigFont{#2pt}}%
  \expandafter\x
    \csname \romannumeral\the\count@ pt\expandafter\endcsname
    \csname @\romannumeral\the\count@ pt\endcsname
  \csname #3\endcsname}%
\fi
\endgroup
\begin{center}
\setlength{\unitlength}{0.00083300in}%
\begingroup\makeatletter\ifx\SetFigFont\undefined
\def\x#1#2#3#4#5#6#7\relax{\def\x{#1#2#3#4#5#6}}%
\expandafter\x\fmtname xxxxxx\relax \def\y{splain}%
\ifx\x\y   
\gdef\SetFigFont#1#2#3{%
  \ifnum #1<17\tiny\else \ifnum #1<20\small\else
  \ifnum #1<24\normalsize\else \ifnum #1<29\large\else
  \ifnum #1<34\Large\else \ifnum #1<41\LARGE\else
     \huge\fi\fi\fi\fi\fi\fi
  \csname #3\endcsname}%
\else
\gdef\SetFigFont#1#2#3{\begingroup
  \count@#1\relax \ifnum 25<\count@\count@25\fi
  \def\x{\endgroup\@setsize\SetFigFont{#2pt}}%
  \expandafter\x
    \csname \romannumeral\the\count@ pt\expandafter\endcsname
    \csname @\romannumeral\the\count@ pt\endcsname
  \csname #3\endcsname}%
\fi
\fi\endgroup
\begin{picture}(3537,3105)(376,-2797)
\thicklines
\put(601,-961){\line( 2, 1){600}}
\put(1201,-661){\line( 2,-1){600}}
\put(2701,-961){\line( 2, 1){600}}
\put(3301,-661){\line( 2,-1){600}}
\put(601,-1411){\line( 3,-1){1192.500}}
\put(1801,-1786){\line( 5, 2){905.172}}
\put(1801,-1411){\line( 5,-2){905.172}}
\put(2701,-1786){\line( 3, 1){1192.500}}
\put(1276,-211){\line( 4, 1){900}}
\put(2176, 14){\line( 4,-1){900}}
\put(1801,-2236){\line( 2,-1){450}}
\put(2251,-2461){\line( 2, 1){450}}
\put(2851,-511){\makebox(0,0)[lb]{\smash{\SetFigFont{12}{14.4}{rm}$(3;\ast,2)$}}}
\put(1051,-511){\makebox(0,0)[lb]{\smash{\SetFigFont{12}{14.4}{rm}$(3;\ast,1)$}}}
\put(376,-1261){\makebox(0,0)[lb]{\smash{\SetFigFont{12}{14.4}{rm}$(2;1,1)$}}}
\put(1576,-1261){\makebox(0,0)[lb]{\smash{\SetFigFont{12}{14.4}{rm}$(2;2,1)$}}}
\put(3676,-1261){\makebox(0,0)[lb]{\smash{\SetFigFont{12}{14.4}{rm}$(2;2,2)$}}}
\put(2476,-1261){\makebox(0,0)[lb]{\smash{\SetFigFont{12}{14.4}{rm}$(2;1,2)$}}}
\put(1576,-2086){\makebox(0,0)[lb]{\smash{\SetFigFont{12}{14.4}{rm}$(1;1,\ast)$}}}
\put(2476,-2086){\makebox(0,0)[lb]{\smash{\SetFigFont{12}{14.4}{rm}$(1;2,\ast)$}}}
\put(1951,164){\makebox(0,0)[lb]{\smash{\SetFigFont{12}{14.4}{rm}$(4;\ast,\ast)$}}}
\put(2026,-2761){\makebox(0,0)[lb]{\smash{\SetFigFont{12}{14.4}{rm}$(0;\ast,\ast)$}}}
\end{picture}
\end{center}
\caption{$P(3,\{[1,2],[2,3]\},2)$}
\label{F_pin1}
\end{figure}
}\end{example}
\begin{example}
\label{Ex_empty}
{\em If $\iv=\emptyset$ then $P(n,\emptyset,N)$ is the same chain
$\0=0<1<2\cdots<n<n+1=\1$ for every positive integer $N$. 
}\end{example}

\begin{proposition}
Let $n$ and $N$ be positive integers and $\iv\df \{I_1,\ldots,I_k\}$ a
nonempty interval system on $[1,n]$. Then $P(n,\iv,N)$ is a graded
poset of rank $n+1$ and we have 
$$f_S(P(n,\iv,N))=N^{|\{j\in[1,k]\::\: S\cap I_j\neq
\emptyset\}|}\quad \mbox{ for all $S\subseteq[1,n]$}.$$
\end{proposition}
\begin{pf}
First we show that $P$ is graded with the rank function given by  
$$\rank((i;p_1,\ldots,p_k))=i.$$ 
Obviously, for every 
$(i;p_1,\ldots,p_k)<(j;q_1,\ldots q_k)$ we have $i<j$.  We only need to
show that whenever $i+1<j$ also holds then there is an element
$(i+1;p'_1,\ldots,p'_k)\in P(n,\iv,N)$ strictly between
$(i;p_1,\ldots,p_k)$ and $(j;q_1,\ldots q_k)$. Let us set  
$$p'_l=\left\{\begin{array}{ll} \ast &\mbox{if $i+1\not\in I_l$,}\\
p_l & \mbox{if $i+1\in I_l$ and $i\in I_l$,}\\
q_l &\mbox{if $i+1\in I_l$ and $j\in I_l$,}\\
\mbox{an arbitrary element of $\{1,2,\ldots,N\}$} & \mbox{if $i+1\in I_l$ and
$i,j\not\in I_l$.}\\
\end{array}\right.$$ 
Observe that whenever we have $\{i,j\}\subseteq I_l$ we also have
$i+1\in I_l$ and $p_l=q_l$, hence there is no contradiction in this
definition of $p'_l$. It is easy to verify that
$(i+1;p'_1,\ldots,p'_k)\in P(n,\iv,N)$ is strictly between
$(i;p_1,\ldots,p_k)$ and $(j;q_1,\ldots q_k)$.

Next we compute $f_S (P(n,\iv,N))$, which is, by definition,
the number of chains $x_1<\cdots <x_{|S|}$
satisfying $\{\rank(x_1),\ldots, \rank (x_{|S|})\}=S$. If for some
$j\in[1,k]$ we have $S\cap I_j=\emptyset$ then every element
$x_t=(i,p_1,\ldots,p_k)$ of such a chain must satisfy $p_j=\ast$. If $S\cap
I_j\neq\emptyset$ then for the elements $x_t=(i,p_1,\ldots,p_k)$
satisfying $\rank (x_t)\in S\cap I_j$ we must have $p_j\in
\{1,2,\ldots,N\}$, and by the definition of the partial order on
$P(n,\iv,N)$ the value of $p_j$ must be the same for all such
$x_t$'s. All other $x_t$'s must satisfy $p_j=\ast$. Conversely, let us fix
a vector $(q_1,\ldots,q_k)$ such that $q_j$ is an arbitrary element of
$\{1,2,\ldots,N\}$ whenever $S\cap I_j\neq \emptyset$ and $q_j=\ast$
otherwise. Then the set $\{ (s;p(s)_1,\ldots,p(s)_k)\::\: s\in S\}$ 
defined by 
$$p(s)_j =\left\{
\begin{array}{ll}
q_j &\mbox{if $s\in I_j$,}\\
 \ast & \mbox{otherwise}\\
\end{array}
\right.
$$
is a chain
in $P(n,\iv,N)$ satisfying $\{ \rank((s;p(s)_1,\ldots,p(s)_k))\::\: s\in
S\}=S$. This shows that $f_S (P(n,\iv,N))$ equals to the number of
possible choices of $(q_1,\ldots,q_k)$, i.e., $N^{|\{j\::\: S\cap
I_j\neq\emptyset\}|}$. 
\end{pf}

\begin{corollary}
\label{C_Pi}
Let $n$ be a positive integer, and $\iv$ an interval system on
$[1,n]$. Then we have  
$$\lim _{N\longrightarrow \infty} \frac{1}{
{f_{[1,n]}(P(n,\iv,N))}}\cdot \sum _{S\subseteq
[1,n]}a_S\cdot f_S(P(n,\iv,N)))=\sum _{{S\in \block_{[1,n]}(\iv)}} a_S.$$
\end{corollary}

In fact, for a non-empty interval system $\iv=\{I_1,\ldots,I_k\}$ we have 
$$ \frac{1}{
{f_{[1,n]}(P(n,\iv,N))}}\cdot \sum _{S\subseteq
[1,n]} a_S\cdot f_S(P(n,\iv,N)))=\frac{1}{
{N^k}}\cdot \sum _{S\subseteq
[1,n]} a_S\cdot N^{|\{j\::\: S\cap
I_j\neq\emptyset\}|}.$$
In this expression $a_S$ is multiplied by $1$ if and only if $S$
intersects every interval of the system $\iv$, otherwise it is
multiplied by a negative power of $N$. Finally, when $\iv=\emptyset$
then, as it was observed in Example \ref{Ex_empty}, every poset
$P(n,\emptyset,N)$ is a chain of rank $n+1$, and we have 
$$\frac{1}{
{f_{[1,n]}(P(n,\emptyset,N))}}\cdot \sum _{S\subseteq
[1,2,\ldots,n]} a_S\cdot f_S(P(n,\emptyset,N)))=\sum _{S\subseteq
[1,n]} a_S$$
immediately.  We can conclude that condition (\ref{E_main}) is necessary:

\begin{corollary}
Suppose $\sum _{S\subseteq [1,n]} a_S\cdot f_S (P)\geq 0$
for every graded poset of rank $n+1$.  Then for every interval system
$\iv$ on $[1,n]$,
$$\sum _{{S\in \block_{[1,n]}(\iv)}} a_S \ge 0.$$
\end{corollary}

Hence we are left to show the sufficiency of (\ref{E_main}).

\begin{proposition}
Assume that the set of coefficients $\{a_S\::\:
S\subseteq [1,n]\}$ satisfies condition (\ref{E_main}). Then we
have 
$$\sum _{S\subseteq [1,n]} a_S\cdot f_S (P)\geq 0$$
for every graded poset of rank $n+1$.
\end{proposition}
\begin{pf}
Let $P$ be a graded poset of rank $n+1$.
For every $i\in[1,n]$ let us fix an arbitrary numbering of
the elements of rank $i$. Given an interval $[p,q]$ of $P$, let
$\phi([p,q])$ denote the first atom in $[p,q]$. (Note that all atoms
of $[p,q]$ have the same rank, namely $\rank(p)+1$.) We will need the following
two elementary observations:
\begin{enumerate}
\item If $y$ covers $x$ then $\phi([x,y])=y$.
\item If $p\in [x,y]\subseteq [x,z]$ and $p=\phi([x,z])$ then $x=\phi([x,y])$. 
\end{enumerate} 
For every $S\subseteq [1,n]$ we define an operation 
$M_S: [1,n]\longrightarrow [1,n]$ by 
$$M_S (i)\df \min \{j\in [i,n+1]\:: j\in S\cup\{n+1\}\}.$$
In other words, $M_S$ assigns to $i$ the smallest element of $S$ which is not
less than $i$, if such an element exists. Otherwise, it assigns $n+1$ to
$i$. 

Consider the set of maximal chains
$$F_S\df\left\{\0=p_0<p_1<\cdots<p_n<p_{n+1}=\1 \::\: \forall
i\in[1,n] \left(p_i=\phi \left([p_{i-1},p_{M_S(i)}]\right)\right)\right\}.
$$
We claim that $F_S$ contains exactly $f_S(P)$ elements. 
For every $i\in S$ we have $M_S(i)=i$ and so, by our second elementary
observation$p_i=\phi\left([p_{i-1},p_{M_S(i)}]\right)
=\phi\left([p_{i-1},p_i]\right)$ is trivially satisfied. Hence
it is sufficient to show 
that every chain $\{p_s\::\: s\in S\}$ satisfying $\{\rank(p_s)\::\:
s\in S\}=S$ may be uniquely extended to a maximal chain
$\{p_0,p_1,\ldots,p_n\}\in F_S$. The only possible choice for $p_0$ is
$\0$. Assume by induction that we have found a unique possible value for
$p_0,\ldots,p_m$. Let $i$ be the smallest rank above $m$ such that
$i\not\in S$. Then the only possible value of $p_i$ is $\phi
([p_{i-1},p_{M_S(i)}])$, and choosing this value we obtain that 
$\{p_1,\ldots,p_i\}\cup \{p_s\::\:s\in S\}$ is a chain. This recursive
algorithm shows that we have at most one extension of $\{p_s\::\: s\in
S\}$ to a maximal chain in $F_S$. On the other hand, at the end of the
algorithm we obtain a maximal chain belonging to $F_S$, since all
defining conditions are satisfied. In particular, we obtain that
$F_{[1,n]}$ contains all maximal chains, while
$F_{\emptyset}$ contains the unique maximal chain for which every
$p_i$  ($i=1,2,\ldots,n$) is the first element among all elements of
rank $i$ covering $p_{i-1}$.  

Let us fix now a maximal chain $C\df
\{\0=p_0<p_1<\cdots<p_n<p_{n+1}=\1\}$, and determine all those sets of
ranks $S\subseteq [1,n]$ for which $C$ belongs to $F_S$. In
view of our two trivial observations, for every $i\in [1,n]$
there is a largest $j\in [i,n+1]$ such that $p_i=\phi
\left([p_{i-1},p_j]\right)$. Let us denote this largest $j$ by $\psi
(C,i)$. Obviously, $p_i=\phi\left([p_{i-1},p_{M_S(i)}]\right)$ is satisfied if
and only if we have $M_S(i)\leq\psi (C,i)$, or equivalently 
$$(S\cup \{n+1\})\cap [i,\psi (C,i)]\neq \emptyset.$$
Introducing 
$$\iv _C\df \{[i,\psi (C,i)]\::\: i\in [1,n],\quad \psi (C,i)\neq n+1\}$$ 
we may say that $C\in F_S$ if an only if $S$ blocks the system $\iv _C$. 

For every $S\subseteq [1,n]$, let us put a weight $a_S$ on
every $C\in F_S$. On the one hand, the sum of all weights put on
the chains is $\sum _{S\subseteq [1,n]} a_S\cdot f_S$. On the
other hand, the total weight associated to an individual chain $C$ is 
$\sum_{S\in \block_{[1,n]}(\iv)} a_S$. Hence we obtain
$$\sum _{S\subseteq [1,n]} a_S\cdot f_S=\sum _{C\in
F_{[1,n]}}\sum_{S\in \block_{[1,n]}(\iv)} a_S,$$  
and so $\sum _{S\subseteq [1,n]} a_S\cdot f_S$ is a sum of
nonnegative terms, if (\ref{E_main}) is satisfied. 
\end{pf}

\section{The Facets of the Cone $\cone{n+1}$}

In this section we show that the inequalities in (\ref{E_main}) give
facets of the cone $\cone{n+1}$. We also show that the
number of facets is a Catalan number.  

\begin{proposition}
\label{P_facet}
Every condition of the form $\sum _{S\in \block_{[1,n]} (\iv)}
a_S\geq 0$ defines a facet of the cone $\cone{n+1}$.
\end{proposition}
\begin{pf}
Assume by way of contradiction that $\sum _{S\in
\block_{[1,n]}(\iv)} a_S\geq 0$ is not a facet.  Then,
by Farkas' lemma, there exists nonnegative numbers $c_1,\ldots,c_k$ and 
interval systems $\iv _1,\ldots,\iv _k$ such that we have
$$\sum _{S\in \block_{[1,n]} (\iv)} a_S=\sum _{i=1}^k c_i\cdot 
\sum _{S\in \block_{[1,n]}(\iv_i)} a_S,$$
and the conditions $\sum _{S\in \block_{[1,n]}(\iv_i)} a_S\geq 0$
are different from $\sum _{S\in \block_{[1,n]}(\iv)} a_S\geq
0$. Since the set $[1,n]$ blocks every interval system, the
coefficient of $a_{[1,n]}$ is $1$ on the left hand side and
$\sum_{i=1}^k c_i$ on the right hand side. Hence we must have 
$$\sum _{i=1} ^k c_i=1.$$
The coefficient of every other $a_S$ is zero or one one the left hand
side, and a convex combination of zeros and ones on the right hand
side. Hence we must have 
$$\block_{[1,n]} (\iv)=\block_{[1,n]} (\iv_1)
= \block_{[1,n]} (\iv_2)=\cdots=\block_{[1,n]} (\iv_k),
$$
contradicting our assumption.
\end{pf}

More directly, Proposition \ref{P_facet} follows from the fact that a
subset of vertices of the $2^n$-dimensional cube are vertices of their
convex hull.

\begin{figure}
\begingroup\makeatletter
\def\x#1#2#3#4#5#6#7\relax{\def\x{#1#2#3#4#5#6}}%
\expandafter\x\fmtname xxxxxx\relax \def\y{splain}%
\ifx\x\y   
\gdef\SetFigFont#1#2#3{%
  \ifnum #1<17\tiny\else \ifnum #1<20\small\else
  \ifnum #1<24\normalsize\else \ifnum #1<29\large\else
  \ifnum #1<34\Large\else \ifnum #1<41\LARGE\else
     \huge\fi\fi\fi\fi\fi\fi
  \csname #3\endcsname}%
\else
\gdef\SetFigFont#1#2#3{\begingroup
  \count@#1\relax \ifnum 25<\count@\count@25\fi
  \def\x{\endgroup\@setsize\SetFigFont{#2pt}}%
  \expandafter\x
    \csname \romannumeral\the\count@ pt\expandafter\endcsname
    \csname @\romannumeral\the\count@ pt\endcsname
  \csname #3\endcsname}%
\fi
\endgroup
\begin{center}
\setlength{\unitlength}{0.00083300in}%
\begingroup\makeatletter\ifx\SetFigFont\undefined
\def\x#1#2#3#4#5#6#7\relax{\def\x{#1#2#3#4#5#6}}%
\expandafter\x\fmtname xxxxxx\relax \def\y{splain}%
\ifx\x\y   
\gdef\SetFigFont#1#2#3{%
  \ifnum #1<17\tiny\else \ifnum #1<20\small\else
  \ifnum #1<24\normalsize\else \ifnum #1<29\large\else
  \ifnum #1<34\Large\else \ifnum #1<41\LARGE\else
     \huge\fi\fi\fi\fi\fi\fi
  \csname #3\endcsname}%
\else
\gdef\SetFigFont#1#2#3{\begingroup
  \count@#1\relax \ifnum 25<\count@\count@25\fi
  \def\x{\endgroup\@setsize\SetFigFont{#2pt}}%
  \expandafter\x
    \csname \romannumeral\the\count@ pt\expandafter\endcsname
    \csname @\romannumeral\the\count@ pt\endcsname
  \csname #3\endcsname}%
\fi
\fi\endgroup
\begin{picture}(2444,2434)(1179,-2173)
\thicklines
\multiput(3601,239)(-8.98876,0.00000){268}{\makebox(6.6667,10.0000){\SetFigFont{10}{12}{rm}.}}
\put(1201,239){\line( 1, 0){2400}}
\put(3601,239){\line( 0,-1){600}}
\put(3601,-361){\line(-1, 0){1200}}
\put(2401,-361){\line( 0,-1){600}}
\put(2401,-961){\line(-1, 0){600}}
\put(1801,-961){\line( 0,-1){600}}
\put(1801,-1561){\line(-1, 0){600}}
\put(1201,-1561){\line( 0, 1){1800}}
\multiput(1201,-1561)(0.00000,-60.00000){11}{\makebox(6.6667,10.0000){\SetFigFont{10}{12}{rm}.}}
\multiput(1201,-2161)(60.00000,0.00000){11}{\makebox(6.6667,10.0000){\SetFigFont{10}{12}{rm}.}}
\multiput(1801,-2161)(0.00000,60.00000){11}{\makebox(6.6667,10.0000){\SetFigFont{10}{12}{rm}.}}
\multiput(1801,-1561)(60.00000,0.00000){11}{\makebox(6.6667,10.0000){\SetFigFont{10}{12}{rm}.}}
\multiput(2401,-1561)(0.00000,60.00000){11}{\makebox(6.6667,10.0000){\SetFigFont{10}{12}{rm}.}}
\multiput(2401,-961)(60.00000,0.00000){11}{\makebox(6.6667,10.0000){\SetFigFont{10}{12}{rm}.}}
\multiput(3001,-961)(0.00000,60.00000){11}{\makebox(6.6667,10.0000){\SetFigFont{10}{12}{rm}.}}
\multiput(1201,-961)(60.00000,0.00000){11}{\makebox(6.6667,10.0000){\SetFigFont{10}{12}{rm}.}}
\multiput(1201,-361)(60.00000,0.00000){21}{\makebox(6.6667,10.0000){\SetFigFont{10}{12}{rm}.}}
\multiput(1801,239)(0.00000,-60.00000){21}{\makebox(6.6667,10.0000){\SetFigFont{10}{12}{rm}.}}
\multiput(2401,239)(0.00000,-60.00000){11}{\makebox(6.6667,10.0000){\SetFigFont{10}{12}{rm}.}}
\multiput(3001,239)(0.00000,-60.00000){11}{\makebox(6.6667,10.0000){\SetFigFont{10}{12}{rm}.}}
\put(1351,-136){\makebox(0,0)[lb]{\smash{\SetFigFont{12}{14.4}{rm}[1,4]}}}
\put(1951,-136){\makebox(0,0)[lb]{\smash{\SetFigFont{12}{14.4}{rm}[2,4]}}}
\put(2551,-136){\makebox(0,0)[lb]{\smash{\SetFigFont{12}{14.4}{rm}[3,4]}}}
\put(3151,-136){\makebox(0,0)[lb]{\smash{\SetFigFont{12}{14.4}{rm}[4,4]}}}
\put(1351,-736){\makebox(0,0)[lb]{\smash{\SetFigFont{12}{14.4}{rm}[1,3]}}}
\put(1351,-1336){\makebox(0,0)[lb]{\smash{\SetFigFont{12}{14.4}{rm}[1,2]}}}
\put(1351,-1936){\makebox(0,0)[lb]{\smash{\SetFigFont{12}{14.4}{rm}[1,1]}}}
\put(1951,-736){\makebox(0,0)[lb]{\smash{\SetFigFont{12}{14.4}{rm}[2,3]}}}
\put(1951,-1336){\makebox(0,0)[lb]{\smash{\SetFigFont{12}{14.4}{rm}[2,2]}}}
\put(2551,-736){\makebox(0,0)[lb]{\smash{\SetFigFont{12}{14.4}{rm}[3,3]}}}
\end{picture}
\end{center}
\caption{Representation of the dual ideal $\{[1,2],[2,3],[4,4]\}^+$ on $[1,4]$.}
\label{F_dual}
\end{figure}
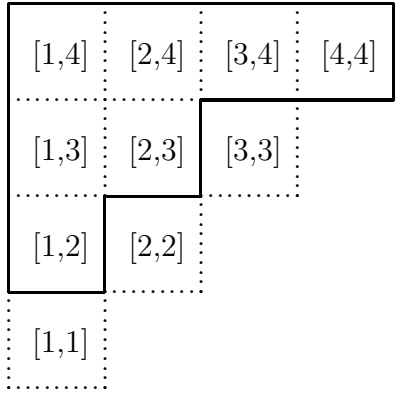

The number of facets is thus equal to the number of families of
sets of the form $\block _{[1,n]}(\iv)$ where $\iv$ is an
arbitrary interval system on $[1,n]$. In order to count the
number of these families, observe that by Proposition \ref{P_dual2}
we may replace $\iv$ by $\iv ^+$. Hence the number of facets of the cone
$\cone{n+1}$ is equal to the number of dual ideals of subsets of $[1,n]$
generated by intervals. There is a bijection between such dual ideals
and Ferrers shapes contained in the shape of the 
partition $(n,n-1,\ldots,1)$, defined as follows. Given an $n\times n$
square, write the interval $[i,j]\subseteq [1,n]$ into the box in the
$j$th row and $i$th column. The boxes into which we have written an
interval form the Ferrers shape of the partition
$(n,n-1,\ldots,1)$. Clearly an interval system $\iv$ on $[1,n]$ is
the family of all intervals in a dual ideal, if and only if for every box
representing an interval $I\in\iv$ all 
boxes above and to the left are marked with an interval from
$\iv$. Equivalently, the boxes representing $\iv$ form a Ferrers shape. Figure
\ref{F_dual} shows the Ferrers shape representation for the dual ideal
$\{[1,2],[2,3],[4,4]\}^+$ on $[1,4]$. 

The number of Ferrers shapes contained in the shape of the
partition $(n,n-1,\ldots,1)$ is equal to the number of those lattice
paths from the lower left corner of the
box marked with $[1,1]$ to the upper right corner of the box marked
with $[n,n]$ which use only unit steps up and to the right and which
never leave the shape of the partition $(n,n-1,\ldots,1)$. 

\begin{corollary}
The number of facets of the cone $\cone{n+1}$ is the
Catalan number\\ $\frac{1}{n+1}\cdot {2\cdot (n+1)\choose n}$. 
\end{corollary} 

Let us also note the following consequence of Corollary \ref{C_block} 
which gives a practically more useful description of the facets of
$\cone{n+1}$   
\begin{corollary}
\label{C_clutter}
Every inequality in (\ref{E_main}) may be written as 
$$\sum _{S\in \block_{[1,n]}(\iv)} a_S\geq 0$$ 
for some antichain of intervals $\iv$ on $[1,n]$.
\end{corollary}

\begin{example}\label{banker}
{\em This inequality was found by Billera and Liu (see
\cite{Billera-Liu}). For every graded poset $P$ of rank $3$  we have
$$f_{1,3}(P)-f_{1}(P)+f_{2}(P)-f_{3}(P)\geq 0.$$
In fact, we may apply Theorem \ref{T_main} for $a_{1,3}=1$, $a_1=-1$,
$a_2=1$, $a_3=-1$, and $a_S=0$ for every other subset $S$ of
$\{1,2,3\}$. There are $14$ dual ideals generated by intervals on the
set $[1,3]$, which are given by the antichain of their minimal intervals
in Table \ref{T_banker}. 
Evaluating $\sum_{S\in\block_{[1,3]}(\iv)} a_S$ for
these $14$ interval systems we see that Condition (\ref{E_main}) is satisfied.
}\end{example}

\begin{table}
\begin{center}
\begin{tabular}{|c|c||c|c|}
\hline
$\iv$ & $\sum_{S\in\block_{[1,3]}(\iv)} a_S$ & $\iv$ & 
$\sum_{S\in\block_{[1,3]}(\iv)} a_S$ \\
\hline
\hline
$\emptyset$ & 0 &  $\{[1],[2,3]\}^+$ & 1 \\
\hline
$\{[1,3]\}^+$ & 0 & $\{[2]\}^+$ & 1 \\
\hline
$\{[1,2]\}^+$ & 1 & $\{[1,2],[3]\}^+$ & 1 \\
\hline
$\{[2,3]\}^+$ & 1 & $\{[1],[2]\}^+$ & 0 \\
\hline
$\{[1]\}^+$ & 0 & $\{[1],[3]\}^+$ & 1 \\
\hline
$\{[1,2],[2,3]\}^+$ &  1 & $\{[2],[3]\}^+$ & 0 \\
\hline
$\{[3]\}^+ $& 0 & $\{[1],[2],[3]\}^+$ & 0 \\
\hline
\end{tabular}
\end{center}
\caption{The value of $\sum_{S\in\block_{[1,3]}(\iv)} a_S$ for
$\sum_{S\subseteq[1,3]} a_S\cdot f_S=f_{1,3}-f_{1}+f_{2}-f_{3}$}
\label{T_banker}
\end{table}
 
\section{Projections and Convolutions}

In this section we present linear projections $\proj
_m^{n+1} :\flag{n+1}\longrightarrow \flag{n}$ which allow us to describe
the cone $\cone{n+1}$ in terms of the cone $\cone{n}$. We show that
these projections are nicely compatible with the convolution operation
defined in \cite{Billera-Liu}.

\begin{definition}
Let $n>0$ and $m\in [0,n]$ be integers. For an arbitrary form\\
$\sum_{S\subseteq [1,n]} a_S\cdot f^{n+1}_S \in\flag{n+1}$ we define its
{\em $m$-th projection into $\flag{n}$} by 
$$\proj^{n+1}_m \left(\sum_{S\subseteq [1,n]} a_S\cdot f^{n+1}_S\right)
\df \sum_{S\subseteq [1,n-1]} \left(\chi ( (S\cup \{0\})\cap [m,n-1]\neq
\emptyset)\cdot a_S + a_{S\cup \{n\}}\right)\cdot f^n_S.$$
We extend this definition to $n=0$ and to negative $m$'s by setting
$\proj^1_0 \left(f^1_{\emptyset}\right)\df 1$
and $\proj ^{n+1}_m\df \proj ^{n+1}_0$ whenever $m<0$. 
\end{definition}
Equivalently, the effect of the projections $\proj^{n+1}_m$ on the
operators $f^{n+1}_S$ may be given by
\begin{equation}
\label{E_projop}
\proj^{n+1}_m \left(f^{n+1}_S\right)=\left\{
\begin{array}{ll}
f^n_{S\setminus\{n\}}&\mbox{if $n\in S$,}\\
\chi ( (S\cup \{0\})\cap [m,n-1]\neq\emptyset)\cdot f^n_S&\mbox{if
$n\not\in S$}.\\
\end{array}
\right.
\end{equation}
In particular, for $m=0$, $(S\cup \{0\})\cap [0,n-1]\neq
\emptyset$ holds for every $S\subseteq [1,n-1]$, and we have  
$$\proj^{n+1}_0 \left(\sum_{S\subseteq [1,n]} a_S\cdot f^{n+1}_S\right)=
\sum_{S\subseteq [1,n-1]} \left(a_S+ a_{S\cup \{n\}}\right)\cdot
f^n_S.$$
and
\begin{equation}
\label{E_projop0}
\proj^{n+1}_0 \left(f^{n+1}_S\right)=f^n_{S\setminus\{n\}}.
\end{equation}
At the other extreme, for $m=n$ the interval $[m,n-1]$ is the empty set,
and we have 
$$\proj^{n+1}_n \left(\sum_{S\subseteq [1,n]} a_S\cdot f^{n+1}_S\right)=
\sum_{S\subseteq [1,n-1]} a_{S\cup \{n\}}\cdot f^n_S.$$
Given an arbitrary interval system $\iv$ on $[1,n-1]$ and $m\in [1,n]$,
a set $T\subseteq [1,n]$ blocks $\iv\cup\{[m,n]\}$ if and only if 
$S\df T\setminus \{n\}\subseteq [1,n-1]$ blocks $\iv$ and either $n\in
T$, i.e., $T=S\cup\{n\}$, or $n\not\in T$ and $S=T$ blocks $[m,n-1]$. 
Moreover, $(S\cup\{0\})\cap [m,n-1]\neq \emptyset$ is equivalent to $S\cap
[m,n-1]\neq\emptyset$.  
Hence we have
\begin{equation}
\label{E_proj}
\sum_{T\in\block _{[1,n]} (\iv \cup \{[m,n]\})} a_T=
\sum_{S\in\block _{[1,n-1]} (\iv)} \left(\chi ( (S\cup \{0\})\cap [m,n-1]\neq
\emptyset)\cdot a_S + a_{S\cup \{n\}}\right).
\end{equation} 
Similarly, for $m=0$, a set $T\subseteq [1,n]$ blocks $\iv$ if and only
if $S\df T\setminus \{n\}\subseteq [1,n-1]$ blocks $\iv$. Moreover,
as noted earlier, $\chi((S\cup\{0\})\cap [0,n-1]\neq \emptyset)=1$ for
every $S\subseteq[1,n-1]$. These observations yield 
\begin{equation}
\label{E_proj0}
\sum_{T\in\block _{[1,n]} (\iv)} a_T=
\sum_{S\in\block _{[1,n-1]} (\iv)} \left(\chi ( (S\cup \{0\})\cap [0,n-1]\neq
\emptyset)\cdot a_S + a_{S\cup \{n\}}\right) .
\end{equation} 
Using equations (\ref{E_proj}) and (\ref{E_proj0}) we may show the
following. 
\begin{theorem}
\label{P_proj}
A form $\sum_{S\subseteq [1,n]}a_S\cdot f^{n+1}_S\in\flag{n+1}$
belongs to $\cone{n+1}$ if and only if for every 
$m\in [0,n]$ the projection $\proj ^{n+1}_m \left( \sum_{S\subseteq [1,n]}a_S\cdot f^{n+1}_S\right)$ belongs to $\cone{n}$.
\end{theorem}
\begin{pf}
The necessity is evident in view of of the equations (\ref{E_proj}) and
(\ref{E_proj0}). To prove sufficiency observe that by Corollary
\ref{C_clutter} it suffices to verify the nonnegativity conditions 
\begin{equation}
\label{E_cl}
\sum_{T\in\block _{[1,n]} ({\mathcal J})} a_T\geq 0
\end{equation}
for every antichain of intervals ${\mathcal J}$ on $[1,n]$. Such an antichain is
either also an interval system on $[1,n-1]$ or the union of the antichain
$\iv\df \{I\in {\mathcal J}\::\: I\subseteq [1,n-1]\}$ and of the singleton
$\{[m.n]\}$ where $[m,n]$ is the unique interval in ${\mathcal J}$
containing $n$. In either case, we may use (\ref{E_proj0}) or
(\ref{E_proj}), respectively, to show that (\ref{E_cl}) holds. 
\end{pf}

We now introduce some linear operators acting on $\flag{}$ which will be
useful in giving a simple expression for $\proj ^{n+1}_m$ for $m>0$.

\begin{definition}
Let $n$ be a positive integer and $k\in [0, n-1]$. We define
$\cut^{n+1}_k:\flag{n+1} \longrightarrow \flag{n}$ by setting 
$$\cut^{n+1}_k (f^{n+1}_S)\df \chi (S\subseteq [1,k])\cdot f^n _{S\cap
[1,k]}$$for every $S\subseteq [1,n]$. We extend this definition to $n=0$
and to negative $k's$ by setting $\cut ^1_0(f^1_{\emptyset})\df 0$ and 
$\rho^{n+1}_k\df 0$ whenever $k<0$.
\end{definition}
In particular, for $n>0$ we have
$$\cut^{n+1}_0 \left(f^{n+1}_S\right)=\left\{\begin{array}{ll}
f^n_{\emptyset}&\mbox{if $S=\emptyset$,}\\
0& \mbox{otherwise.}\\
\end{array}
\right.
$$
It is an easy consequence of (\ref{E_projop}) that we have
\begin{equation}
\label{E_pc}
\proj^{n+1}_m=\proj ^{n+1}_0 -\cut^{n+1}_{m-1}\quad\mbox{for $n\in \Nnn$
and $m\in [0,n]$.}
\end{equation}

Let us consider now the effect of the projection operations
$\proj^{n+1}_m$ on a convolution of two chain operators. 
\begin{proposition}
\label{P_prod0}
We have 
$$\proj ^{m+n}_0 \left(f^m_S\ast f^n_T\right)=f^m_S\ast \proj^n_0
\left(f^n_T\right) 
$$
for all positive $m,n$ and sets $S\subseteq [1,m-1]$, $T\subseteq
[1,n-1]$.
\end{proposition}
\begin{pf}
Assume first that we have $n\geq 2$. A simple substitution into the
definitions and (\ref{E_projop0}) yields 
\begin{eqnarray*}
\proj ^{m+n}_0 \left(f^m_S\ast f^n_T\right)
&=& \proj ^{m+n}_0 \left(f^{m+n}_{S\cup\{m\}\cup(T+m)}\right)
= f^{m+n-1}_{S\cup\{m\}\cup(T+m)\setminus\{m+n-1\}}\\
&=& f^{m+n-1}_{S\cup\{m\}\cup(T\setminus\{n-1\})+m}
=f^m_S\ast \proj^n_0\left(f^n_T\right).
\end{eqnarray*}
For $n=1$ we have
$$
\proj ^{m+1}_0 \left(f^m_S\ast f^1_{\emptyset}\right)
= \proj ^{m+1}_0 \left(f^{m+1}_{S\cup\{m\}}\right)
= f^{m}_{S\cup\{m\}\setminus\{m\}}\\
= f^{m}_{S}=f^{m}_{S}\ast 1
=f^m_S\ast \proj^n_0\left(f^1_{\emptyset}\right)
$$
by $\proj^n_0\left(f^1_{\emptyset}\right)=1$.
\end{pf}

\begin{proposition}
\label{P_prodcut}
Let $m$ and $n$ be positive integers, $S\subseteq [1,m-1]$, $T\subseteq
[1,n-1]$. Then for every $k\in [0,m+n-2]$ we have
$$\cut^{m+n}_k\left(f^m_S\ast f^n_T\right)=f^m_S\ast
\cut^n_{k-m}\left(f^n_T\right).$$  
\end{proposition}
\begin{pf}
If $k$ is negative then both sides are identically zero. Assume $0\leq
k\leq m-1$. Then we have
\begin{eqnarray*}
\cut^{m+n} _k\left(f^m_S\ast f^n_T\right)
&=&\cut^{m+n}_k\left(f^{m+n}_{S\cup\{m\}\cup(T+m)}\right)\\
&=& \chi (S\cup\{m\}\cup(T+m)\subseteq [1,k])\cdot
f^{m+n-1}_{\left(S\cup\{m\}\cup(T+m)\right)\cap [1,k]}=0
\end{eqnarray*}
since $m\not\in [1,k]$. 
On the other hand, $\rho^n_{k-m}$ is identically zero by definition, and
so the equality of the both sides holds trivially.

Assume finally that $m\leq k\leq m+n-2$ holds. Then we must have $n\geq
2$, $S$ is a subset of $[1,k]$ and $S\cup\{m\}\cup(T+m)\subseteq
[1,k]$ is equivalent to $T\subseteq [1,k-m]$. Therefore we have
\begin{eqnarray*}
\cut^{m+n}_k\left(f^m_S\ast f^n_T\right)
&=&\cut^{m+n}_k\left(f^{m+n}_{S\cup\{m\}\cup(T+m)}\right)\\
&=&\chi (S\cup\{m\}\cup(T+m)\subseteq [1,k])\cdot
f^{m+n-1}_{\left(S\cup\{m\}\cup(T+m)\right)\cap [1,k]}\\
&=&\chi (T\subseteq [1,k-m])\cdot f^{m+n-1}_{S\cup\{m\}\cup
(T\cap[1,k-m]+m)}\\ 
&=& \chi (T\subseteq [1,k-m])\cdot f^m_S\ast f^{n-1}_{T\cap [1,k-m]}\\
&=&f^m_S\ast \cut^n_{k-m} \left(f^n_T\right).
\end{eqnarray*}
\end{pf}

As a corollary of Propositions \ref{P_prod0} and \ref{P_prodcut}, and of
equation (\ref{E_pc}) we obtain 
\begin{corollary}
\label{C_prod} The equality
$$
\proj^{m+n}_k \left(f^m_S\ast f^n_T\right)=f^m_S\ast
\proj^{n}_{k-m}\left(f^n_T \right)$$
holds for $m,n>0$, $k\in [0,m+n-1]$, $S\subseteq [1,m-1]$, and
$T\subseteq [1,n-1]$. 
\end{corollary} 

In section \ref{S_extreme} we will need to consider
the maximum element of the {\em support} of
a form $\sum_{S\subseteq [1,n]} a_S\cdot f^{n+1}_S$. 

\begin{definition}
The {\em support} of a form $F=\sum_{S\subseteq [1,n]} a_S\cdot
f^{n+1}_S\in\flag{n+1}$ is the family of sets 
$$\supp (F)\df \{S\subseteq [1,n]\::\: a_S\neq 0\}.$$
We call the maximum element of $\bigcup _{S\in \supp (F)} S$ the {\em
largest letter occurring in $F$}. 
\end{definition} 

\begin{lemma}
The projections $\proj^{n+1}_m$ do not increase the largest letter occurring
in a form.
\end{lemma}
\begin{pf}
Let us denote the largest occurring letter in 
$F= \sum_{S\subseteq [1,n]} a_S\cdot f^{n+1}_S\in\flag{n+1}$ by
$l$. If $l=n$ then there is nothing to prove, so we may assume $l<n$.
Evidently $k\in [1,n]$ is greater than $l$ if and
only if for every $S\subseteq [1,n]$ containing $k$ we have $a_S=0$. 
Hence it is sufficient to show the following: if for every $S\subseteq
[1,n]$ containing some $k>l$ we have $a_S=0$ then for every $S\subseteq
[1,n-1]$ containing some $k>l$ we have 
$$\chi ( (S\cup \{0\})\cap [m,n-1]\neq
\emptyset)\cdot a_S + a_{S\cup \{n\}}=0.$$
The second term in this sum is zero since $S\cup\{n\}$ contains the
letter $n$ which larger than $l$ by assumption. The first term is zero too,
since $S$ contains some $k>l$. 
\end{pf}

\section{Extreme Rays of $\cone{n+1}$}
\label{S_extreme}

In this section we study the extreme rays of the cones $\cone{n+1}$ where $n\in
\Nnn$. First we show that $f^{n+1}_{\emptyset}$ is an extreme ray of the cone
$\cone{n+1}$ for every $n\geq 0$. Then we show for $k=1,2,\ldots,n$ that
the inclusions 
\begin{equation}
\label{E_shiftop}
\begin{array}{rrcl}
\shift{k}^{n+1}:&\flag{n+1}&\longrightarrow&\flag{n+2}\\
&f^{n+1}_S&\longmapsto & f^{n+2}_{\shift{k}(S)}\\ 
\end{array}
\end{equation}
induced by the {\em shift operators}  
\begin{equation}
\label{E_shift}
\begin{array}{rrcl}
\shift{k}:&\Ppp&\longrightarrow&\Ppp\\
&i&\longmapsto & 
\left\{\begin{array}{ll}
i&\mbox{if $i \leq k$}\\
i+1&\mbox{if $i\geq k+1$}\\
\end{array}
\right.
\\
\end{array}
\end{equation}
embed $\cone{n+1}$ into $\flag{n+2}$ as a face of
$\cone{n+2}$. Finally we
describe completely those situations when the convolution $F\ast G$ of
an extreme ray $F$ of the cone $\cone{m}$ and an extreme ray $G$ of the
cone $\cone{n}$ is an extreme ray of the cone $\cone{n+m}$.  The main
result of this section is the following theorem.  
\begin{theorem}
\label{T_per}
Let $F\in\cone{m}$ and $G\in\cone{n}$ be extreme rays in their
respective cones. Then the convolution $F\ast G\in\cone{m+n}$ is an
extreme ray, unless $F=F'\ast f^k_{\emptyset}$ and
$G=f^l_{\emptyset}\ast G'$ for some $k \le m$, $l\le n$, $F'\in
\cone{m-k}$, and $G'\in \cone{n-l}$.
\end{theorem}

To interpret this result, note that Proposition \ref{P_conv} is
equivalent to the following.
\begin{proposition}
\label{P_convc}
Let $F\in \flag{m}$ and $G\in\flag{n}$. The convolution $F\ast
G\in\flag{m+n}$ belongs to $\cone{n+m}$ if and only if exactly one of
the following holds:
\begin{enumerate}
\item[(i)] $F\in\cone{m}$ and $G\in\cone{n}$, or 
\item[(ii)] $-F\in\cone{m}$ and $-G\in\cone{n}$.
\end{enumerate}
\end{proposition}
Thus the convolution of two extreme rays is surely in the cone of valid
inequalities. If, say, $G$ is a positive linear combination of
$G'$ and $G''$ from $\cone{n}$ then $F\ast G$ will be the positive linear
combination of $F\ast G'$ and $F\ast G''$ from $\cone{m+n}$. (Recall
that according to Proposition \ref{P_free} the ring
$A$ is a free associative algebra, and so it has no zero divisors.) Thus only
the convolution of extreme rays may yield an extreme ray. It may happen
that the convolution of extremes is not extreme: for every $m,n>0$ the
operators $f^m_{\emptyset}$ and
$f^n_{\emptyset}$ are extreme rays in $\cone{m}$ and $\cone{n}$ respectively,
yet
$$f^m_{\emptyset}\ast
f^n_{\emptyset}=f^{m+n}_m
=\left(f^{m+n}_m-f^{m+n}_{\emptyset}\right)+f^{m+n}_{\emptyset}$$
where both $f^{m+n}_m-f^{m+n}_{\emptyset}$ and
$f^{m+n}_{\emptyset}$ belong to $\cone{m+n}$. Theorem \ref{T_per} affirms
that, essentially, only such anomalies may occur. 

According to Proposition \ref{P_ufact} the semigroup of homogeneous
polynomials of $\flag{}$ has unique factorization. In view of Proposition
\ref{P_convc}, an expression $F\in\flag{n}$ belongs to $\cone{n}$ if and
only if every factor $F_i\in \flag{n_i}$ in its complete
homogeneous factorization $F=F_{1}\ast\cdots\ast F_{k}$ may be chosen to
belong to $\cone{n_i}$. (Since uniqueness holds up to a choice of
nonzero constant factors, we may change the signs of the factors of two
$F_i$'s at a time.) When $F$ is an extreme ray of $\cone{n}$ then every
$F_i$ must be an extreme ray in its cone, and no two consecutive factors
can be of the form $F_i=f^{n_i}_{\emptyset}$,
$F_{i+1}=f^{n_{i+1}}_{\emptyset}$. Theorem \ref{T_per} implies that the
converse is true as well: every convolution of extreme rays such that no
two consecutive factors are of the form $F_i=f^{n_i}_{\emptyset}$,
$F_{i+1}=f^{n_{i+1}}_{\emptyset}$, is an extreme ray.

One way of showing $F\in\cone{n+1}$ is an
extreme ray is to give $2^n-1$ linearly independent facets of
$\cone{n+1}$ containing
$F$. (Note that since $\cone{n+1}$ contains the positive orthant of
$\flag{n+1}$, it is full dimensional.) Hence it is useful to introduce
the following operators.

\begin{definition}
\label{D_eval}
Let $P$ be an arbitrary graded poset of rank $n+1$, and $\iv$ an
arbitrary interval system on $[1,n]$.  
The operators $\eval{P}^{n+1},\eval{\iv}^{n+1}\in\dflag{n+1}$ 
are given by
$$\eval{P}^{n+1}(f^{n+1}_S)\df f^{n+1}_S(P),\quad\mbox{and}$$ 
$$\eval{\iv}^{n+1}(f^{n+1}_S)\df 
\left\{
\begin{array}{ll}
1&\mbox{if $\forall I\in\iv ~(S\cap I\neq \emptyset)$,}\\
0&\mbox{otherwise.}\\
\end{array}
\right.$$
\end{definition}

As a consequence of this definition we have for $F=\sum
_{S\subseteq [1,n]} a_S\cdot f^{n+1}_S\in\flag{n+1}$ that 
$$\eval{\iv}^{n+1}(F)=\sum_{S\in\block{[1,n]}(\iv)} a_S,$$
hence the hyperplanes determining the facets of $\cone{n+1}$ are
the kernels of the $\eval{\iv}^{n+1}$'s. Let us also note that     
Corollary \ref{C_Pi} may be rewritten as 
$$\eval{\iv}^{n+1}=\lim _{N\longrightarrow \infty} \frac{1}{
{f_{[1,n]}(P(n,\iv,N))}}\cdot \eval{P(n,\iv,N)}^{n+1}.$$

It is easy to show that the operators $\eval{\iv}^{n+1}$, where $\iv$ runs
over all interval systems on $[1,n]$, contain a basis of the vector space
$\dflag{n+1}$. 

\begin{definition}
For every $S\subseteq [1,n]$ let $\iv_S$ denote the interval
system $\{ \{s\}\::\: s\in S\}$. As a shorthand for $\eval{\iv_S}^{n+1}$ we
will use $\eval{S}^{n+1}$.
\end{definition}

\begin{lemma}
\label{L_is}
For every $S,T\subseteq [1,n]$ we have 
$$\eval{S}^{n+1}\left(f^{n+1}_T\right)=\left\{
\begin{array}{ll}
1 &\mbox{if $S\subseteq T$,}\\
0&\mbox{otherwise}.\end{array}\\
\right.$$
\end{lemma}
The proof is straightforward.
\begin{corollary}
\label{C_indep}
The set $\{\eval{S}^{n+1}\::\: S\subseteq [1,n]\}$ is a basis of
$\dflag{n+1}$.  
\end{corollary}

\begin{proposition}
The chain operator $f^{n+1}_{\emptyset}$ is an extreme ray of the cone
$\cone{n+1}$ for every $n\geq 0$.
\end{proposition}
\begin{pf}
Since $f^{n+1}_{\emptyset} (P)>0$ holds for every partially ordered set 
of rank $n+1$, we have $f^{n+1}_{\emptyset}\in \cone{n+1}$. 
We only need to show that $f^{n+1}_{\emptyset}$ lies on at least $2^n-1$
linearly independent facets of $\cone{n+1}$. This is true, since by
Lemma \ref{L_is} we have $\eval{S}^{n+1} (f^{n+1}_{\emptyset})=0$, whenever
$S$ is not the empty set.
\end{pf}

The facets of the form $\KER (\eval{S}^{n+1})$ are also useful in proving
the following proposition.

\begin{proposition}
Let $n>0$ be an integer and $k\in [1,n]$. Then the set
$\shift{k}^{n+1}(\cone{n+1})$ is a face of $\cone{n+2}$.
\end{proposition}
\begin{pf}
Evidently we have 
$$\IM (\shift{k}^{n+1})=\bigcap _{{S\subseteq
[1,n+1]\atop k\in S}} \KER (\eval{S}^{n+1}),$$
since both vector spaces are spanned by those chain operators
$f^{n+2}_S$ for which $k\not\in S$. Thus we only need to show that
$\shift{k}^{n+1}(\cone{n+1})\subseteq \cone{n+2}$. The cone
$\shift{k}^{n+1}(\cone{n+1})$ is then contained in the $2^n$-dimensional face 
$$ \cone{n+2}\cap \left(\bigcap _{{S\subseteq
[1,n+1]\atop k\in S}} \KER (\eval{S}^{n+1})\right)$$ 
of $\cone{n+2}$, and, having
the same dimension, it is also equal to it. Equivalently, we have
to prove that 
$$\eval{\iv}^{n+2}(\shift{k}^{n+1}(F))\geq0$$ 
holds for every interval system $\iv$ on $[1,n+1]$.

Assume first that $\iv$ contains the interval $\{k\}$. Then we have 
$$\eval{\iv}^{n+2}(\shift{k}^{n+1}(F))=0$$
for every $F\in\flag{n+1}$. If $\iv$ does not contain $\{k\}$ then
consider the interval system
$$\iv'\df\{(I\cap [1,k-1])\cup (I\cap [k,n+1]-1)\::\: I\in \iv\}$$
on $[1,n]$. It is easy to verify that we have 
$$\eval{\iv}^{n+2}(\shift{k}^{n+1}(F))=\eval{\iv'}^{n+1}(F)$$
for every $F\in \flag{n+1}$, and so $F\in\cone{n+1}$ implies
$\eval{\iv}^{n+2}(\shift{k}^{n+1}(F))\geq 0$. 
\end{pf}

\begin{corollary}
\label{C_lift}
Given $n\geq 1$ and $k\in [1,n]$, the form $F\in\flag{n+1}$ belongs 
to $\cone{n+1}$ if and only if $\shift{k}^{n+1}(F)$ belongs to
$\cone{n+2}$. Moreover, $F\in\cone{n+1}$ is an extreme ray if and only
if $\shift{k}^{n+1}(F)\in \cone{n+2}$ is an extreme ray.
\end{corollary}

Corollary \ref{C_lift} implies that every extreme ray $F\in\cone{n+2}$
with $\supp(F)\neq [1,n+1]$ is obtained by lifting an extreme ray of 
$\cone{n+1}$ using an embedding $\shift{k}^{n+1}$. Iterated use of
Corollary \ref{C_lift} yields the following.

\begin{theorem}
\label{C_select}
Let $F=\sum _{S\subseteq [1,n]}a_s\cdot f^{n+1}_S \in \flag{n+1}$ be a
form with support $\supp (F)=\{i_1,i_2,\ldots,i_k\}$ and let
$\gamma :[1,k]\longrightarrow \supp (F)$ be the bijection $j\mapsto i_j$. Then 
$F$ belongs to $\cone{n+1}$ if and only if 
$$F'\df \sum_{S\subseteq \supp(F)} a_{\gamma^{-1}(S)}\cdot
f^{k+1}_{\gamma^{-1}(S)}$$ 
belongs to $\cone{k+1}$. Moreover $F\in \cone{n+1}$ is an extreme ray if
and only if $F'\in \cone{k+1}$ is an extreme ray.
\end{theorem}
It follows from Theorem \ref{C_select} and our calculation of the rank
2 case in subsection \ref{ss_gpo} that the functionals $h_i \df f_i -
f_\emptyset$, $1\le i\le n$, will all be extreme in rank $n+1$.

The crucial step in the proof of Theorem \ref{T_per} is the following. 
\begin{lemma}
\label{L_sep}
Let $F=\sum_{S\subseteq [1,n]} a_S\cdot f^{n+1}_S\neq
f^{n+1}_{\emptyset}$ be an extreme ray of 
$\cone{n+1}$ and let $m$ be the largest letter occurring in $F$. Assume
that $F$ cannot be written as $G\ast f^{n+1-m}_{\emptyset}$ for some
$G\in K_m$. Then
there exists integers $k,l$ satisfying $0\leq k<l \leq m$ and an
interval system $\iv$ on $[1,m-1]$ such that  
$$\sum_{S\cup\{ 0\}\in\block{[0,n]}(\iv\cup\{[k,m]\}} a_S=0 \quad
\mbox{and}\quad
\sum_{S\in\block{[1,n]}(\iv\cup\{[l,m]\}} a_S\neq 0\quad\mbox{hold.}$$
\end{lemma}
\begin{pf}
Let us show first that, without loss of generality, we may restrict
ourselves to the case $m=n$. Since $m$ is the largest occurring letter, 
we have 
$$F=\sum _{S\subseteq [1,n]} a_S\cdot f^{n+1}_S
=\sum _{S\subseteq [1,m]} a_S\cdot f^{n+1}_S.$$
The form $F$ is of the form $G\ast f^{n+1-m}_{\emptyset}$ for some
$G\in \cone{m}$ if and only if the form 
$$\widetilde{F}\df \sum _{S\subseteq [1,m]} a_S\cdot f^{m+1}_S$$ 
is of the form $G\ast f^1_{\emptyset}$ for the same $G\in \cone{m}$. 
Given an arbitrary interval system ${\mathcal J}$ on $[1,m]$ we have 
$$\eval{\mathcal J} ^{n+1}(F)=\eval{\mathcal J} ^{m+1}(\widetilde{F}).$$
Therefore, if the lemma holds for the case $m=n$ then we may use the result on
$\widetilde{F}$  to prove the same result for $F$. 

Hence we may assume that $m=n$ holds. The statement is equivalent
to saying that there exists a facet $\KER\left(\eval{\iv}\right)$ of
$\cone{n}$ and integers $0\leq k<l\leq n$ such that $\proj ^{n+1}_k (F)$
belongs to this facet but  $\proj ^{n+1}_l (F)$ does not. 

In the contrary event every facet of $\cone{n}$ containing some $\proj
^{n+1}_k (F)$ also contains all forms $\proj ^{n+1}_l (F)$ for every
$l>k$. Let $G\in\flag{n}$ be an extreme ray of the intersection of all
facets containing $\proj ^{n+1}_n (F)$. (Note that $n$ being the largest
occurring letter, $\proj^{n+1}_n(F)=\sum _{S\subseteq [1,n-1]}
a_{S\cup\{n\}} \cdot f^n_S$ is not the zero form.) By our assumptions,
this form is also an extreme ray of the intersection of all facets
containing $\proj ^{n+1}_i (F)$ for $i=0,1,\ldots,n-1$. Hence we may
choose a small positive number $q\in \Rrr$ such that 
$$ \proj ^{n+1}_i (F)-q\cdot G\in \cone{n}\quad\mbox{holds for
$i=0,1,\ldots,n$}.$$ 
By Corollary \ref{C_prod} we have $\proj^{n+1}_i (G\ast
f^1_{\emptyset})=G$ for $i=1,2,\ldots,n$ and so we obtain 
$$ \proj ^{n+1}_i \left(F-q\cdot G\ast f^1_{\emptyset})\right)\in
\cone{n}$$
for $i=0,1,\ldots,n$. By Theorem \ref{P_proj} we obtain that
$F-q\cdot G\ast f^1_{\emptyset}$ belongs to $\cone{n+1}$. Since $F$ is
an extreme ray, and by Proposition \ref{P_conv} we have $G\ast
f^1_{\emptyset}\in \cone{n+1}$, $F$ must be equal to a nonzero constant
multiple of  $G\ast f^1_{\emptyset}$, contrary to our assumptions. 
\end{pf}
\begin{corollary}
\label{C_sep}
Let $F\in \cone{n+1}$ be an extreme satisfying the conditions of Lemma 
\ref{L_sep}. Then there exists integers $k,l$ satisfying $0\leq k<l \leq
m$ and an interval system $\iv$ on $[1,m-1]$ such that  
$$\sum_{S\cup \{ 0\}\in\block{[0,n]}(\iv\cup\{[k,n]\}} a_S=0 \quad
\mbox{and}\quad
\sum_{S\in\block{[1,n]}(\iv\cup\{[l,n]\}} a_S\neq 0\quad\mbox{hold.}$$
\end{corollary}
In fact, replacing the intervals $[k,m]$ and $[l,m]$ with $[k,n]$ and
$[l,n]$ respectively does not change the sums involved, since whenever
$S$ contains a letter larger than $m$, $a_S$ is zero.

We now proceed to the proof of Theorem \ref{T_per}. Since applying the
chain operators to the dual of every poset yields an anti-isomorphism of
the graded ring $A$ which sends products of the form $F\ast f^l_{\emptyset}$
into products of the form $f^l_{\emptyset}\ast \widehat{F}$, it is
sufficient to show the following ``half'' of the original statement.
\begin{proposition}
Assume $F\in \cone{m}$ and $G\in\cone{n}$ are extreme rays such that 
$F$ is not a convolution of the form $F'\ast f^l_{\emptyset}$. Then
$F\ast G$ is an extreme ray of $\cone{m+n}$.
\end{proposition}
\begin{pf}
By Proposition \ref{P_convc} we know that $F\ast G$ belongs to
$\cone{m+n}$. In order to show that it is an extreme ray, it is
sufficient to find $2^{m+n-1}-1$ interval systems $\iv_1,\iv_2,\ldots,
\iv_{2^{m+n-1}-1}$ such that the operators
$\eval{\iv_1}^{m+n},\eval{\iv_2}^{m+n}, \ldots,
\eval{\iv_{2^{m+n-1}-1}}^{m+n}$ are linearly independent, and they all
vanish on $F\ast G$.  

Assume we have 
$$F=\sum_{S\subseteq [1,m-1]} a^m_S\cdot f^m_S 
\quad\mbox{and}\quad
G=\sum_{S\subseteq [1,n-1]} a^n_S\cdot f^n_S.$$
Then 
$$F*G=\sum _{S\subseteq [1,m+n-1]} a^{m+n}_S\cdot f^{m+n}_S$$
is given by the following formula:
$$
a_S^{m+n}=\left\{
\begin{array}{ll}
a^m_{S\cap [1,m-1]}\cdot a^n_{(S\cap [m+1,m+n-1])-m }&\mbox{if $m\in S$,}\\
0&\mbox{if $m\not\in S$.}\\
\end{array}
\right.
$$
Since for every nonzero $a_S^{m+n}$ we must have $m\in S$, for an
arbitrary interval system $\iv$ on $[1,m+n-1]$ we have 
$$\eval{\iv}^{m+n}(F\ast G)=\eval{\{I\in\iv\::\: m\not\in
I\}}^{m+n}(F\ast G),$$
i.e., we may remove every interval containing $m$ from $\iv$ without
changing the effect of $\eval{\iv}^{m+n}$ on $F\ast G$. The remaining
intervals are either contained in $[1,m-1]$ or in
$[m+1,m+n-1]$. Introducing 
$$\iv'\df \{I\in\iv\::\: I\subseteq [1,m-1]\}\quad\mbox{and}\quad
\iv''\df \{I-m\::\:I\in\iv, I\subseteq [m+1,m+n-1]\}$$ 
we obtain
\begin{equation}
\label{E_eprod}
\eval{\iv}^{m+n}(F\ast G)=\eval{\iv'}^m (F)\cdot \eval{\iv''}^n(G).
\end{equation}

Since $F\in\cone{m}$ is an extreme ray, there exist interval
systems $\iv'_1,\iv'_2,\ldots,\iv'_{2^{m-1}}$ on $[1,m-1]$ such that 
the operators $\eval{\iv'_1}^m,\eval{\iv'_2}^m,\ldots,\eval{\iv'_{2^{m-1}}}^m$ 
are linearly independent, $\eval{\iv'_j}^m (F)=0$ for $j\leq
2^{m-1}-1$, and $\eval{\iv'_{2^{m-1}}}^m (F)$ is strictly positive. 
Similarly, the fact of $G\in\cone{n}$ being an extreme ray implies that 
there exist interval systems $\iv''_1,\iv''_2,\ldots,\iv''_{2^{n-1}}$ on
$[1,n-1]$ such that the operators
$\eval{\iv''_1}^n,\eval{\iv''_2}^n,\ldots,\eval{\iv''_{2^{n-1}}}^n$ are
linearly independent, $\eval{\iv''_j}^n (F)=0$ holds for $j\leq
2^{n-1}-1$, and $\eval{\iv'_{2^{n-1}}}^n (F)$ is strictly positive.
Moreover, since $F$ is not of the form $F'\ast f^l$, by Corollary
\ref{C_sep} we may assume that $\iv'_{2^{m-1}}$ is of the form ${\mathcal
J}'\cup\{[l,m-1]\}$ where ${\mathcal J}'$ is an interval system on
$[1,m-1]$, $l\in [1,m-1]$ and either we have $\eval{{\mathcal J}'}^m(F)=0$,
or there exist a $k\in [1,l-1]$ such that $\eval{{\mathcal J}'\cup
\{[k,m-1]\}}^m(F)=0$ holds. 

Consider now the following interval systems:
\begin{enumerate}
\item[(i)] All interval systems of the form $\iv'_i\cup (\iv''_j+m)$,
where at least one of $i\neq 2^{m-1}$ and $j\neq 2^{n-1}$ holds. 
\item[(ii)]  All interval systems of the form $\iv'_i\cup\{m\}\cup
(\iv''_j+m)$, where at least one of $i\neq 2^{m-1}$ and $j\neq 2^{n-1}$ holds. 
\item[(iii)] The interval system ${\mathcal J}'\cup \{[l,m]\} \cup
(\iv''_{2^{n-1}}+n)$ if we have  $\eval{{\mathcal J}'}^m(F)=0$, or the interval
system ${\mathcal J}'\cup \{[k,m-1]\}\cup \{[l,m]\} \cup
(\iv''_{2^{n-1}}+m)$ if $\eval{{\mathcal J}'\cup \{[k,m-1]\}}^m(F)=0$ holds.
\end{enumerate}
The above list contains $2^{m+n-1}-1$ interval systems: $2^{m+n-2}-1$ of
them are of type (i), $2^{m+n-2}-1$ of them are of type (ii), and there is
exactly one system listed at item (iii). It is easy to see using
(\ref{E_eprod}) that they all vanish on $F\ast G$. We are left to show
that the operators they define are linearly independent. 

First we show that the operator defined by the last system is not in the
span of the operators defined by all others. For this purpose consider
the form $H\df \sum _{S\subseteq [1,m+n-1]} \widetilde{a}_S\cdot
f^{m+n}_S$ given by 
$$\widetilde{a}_S=
\left\{
\begin{array}{ll}
0&\mbox{if $m\in S$,}\\
a_{S\cap [1,m-1]}^m\cdot a_{(S\cap [m+1,m+n-1])-m}^n 
& \mbox{if $m\not\in S$.}\\ 
\end{array}
\right.
$$
It is easy to see that we have 
$$\eval{\iv'_i\cup (\iv''_j+m)}^{m+n}(H)=\eval{\iv'_i}^m (F)\cdot
\eval{\iv''_j}^n (G) $$ 
and so the operators defined by the systems of type (i) vanish on
$H$. The operators defined by the systems of type (ii)
vanish on $H$ too, since we have $\widetilde{a}_S=0$
whenever $m\in S$. By the same reason we also have 
$$\eval{{\mathcal J}'\cup \{[k,m-1]\}\cup \{[l,m]\} \cup
(\iv''_{2^{n-1}}+m)}^{m+n}(H)
=\eval{{\mathcal J}'\cup \{[k,m-1]\}\cup \{[l,m-1]\} \cup
(\iv''_{2^{n-1}}+m)}^{m+n}(H)$$
Since $k$ is less than $l$, the dual ideal of intervals generated by 
${\mathcal J}'\cup \{[k,m-1]\}\cup \{[l,m-1]\} \cup (\iv''_{2^{n-1}}+m)$ 
is the same as the dual ideal generated by 
${\mathcal J}'\cup \{[l,m-1]\} \cup (\iv''_{2^{n-1}}+m)=\iv'_{2^{m-1}}\cup
(\iv''_{2^{n-1}}+m)$. Hence we have 
$$\eval{{\mathcal J}'\cup \{[k,m-1]\}\cup \{[l,m]\} \cup
(\iv''_{2^{n-1}}+m)}^{m+n}(H)= \eval{\iv'_{2^{m-1}}}^m (F)\cdot
\eval{\iv''_{2^{n-1}}}^m (G)\neq 0.$$
Next we show that the intersection of the subspace generated by the
operators associated to the interval systems listed in (i) with the subspace 
generated by the operators associated to the interval systems listed in
(ii) is zero. For this purpose observe that $\flag{n}$ may be written as
a direct sum of two $2^{m+n-2}$-dimensional vector spaces, 
$\flag{n} = A'_n \oplus A''_n$, where
$$A'_n \df \langle f^{m+n}_{S\cup\{m\}}-f^{m+n}_S\::\: S\subseteq [1,m+n-1],\,
m\not\in S\rangle$$
and
$$A''_n \df \langle f^{m+n}_{S\cup\{m\}}\::\: S\subseteq [1,m+n-1],\,
m\not\in S\rangle.$$
Evidently, linear combinations of operators
associated to interval systems of type (i) vanish on $A'_n$
while linear combinations of operators associated to
interval systems of type (ii) have the same values on the
respective generators of $A'_n$ and $A''_n$.
Assume now that $\varepsilon\in \dflag{n}$ is
simultaneously a linear combination of operators associated to
interval systems of type (i) and of operators associated to
interval systems of type (ii). Then $\varepsilon$ vanishes
on $A'_n$ and hence on $A''_n$, so  $\varepsilon=0$.

Finally, let $F_1,\dots,F_{2^{m-1}} \in \flag{m}$ and
$G_1,\dots,G_{2^{n-1}} \in \flag{n}$ be dual bases to
$\{\eval{\iv'_i}^m\}$ and $\{\eval{\iv''_i}^n\}$, respectively,
that is, $\eval{\iv'_i}^m(F_j) = \delta_{ij}$ and 
$\eval{\iv''_k}^m(G_l) = \delta_{kl}$.
Then by (\ref{E_eprod}) we have
$$\eval{\iv'_i\cup (\iv''_k+m)}(F_j \ast G_l) =
\eval{\iv'_i\cup \{m\}\cup (\iv''_k+m)}(F_j \ast G_l) =
\delta_{ij}\delta_{kl},$$
showing the interval systems of types (i) and (ii) to be linearly
independent.
\end{pf}

We can describe all extremes of the cone $\cone{n}$ for $n \le 5$.
Equivalently, these represent the strongest linear inequalities holding
for the flag $f$-vectors of graded posets of these ranks.

\begin{description}
\item[$n=1$] The only extreme ray of $\cone{1}$ is 
$h^1_\emptyset = f^1_\emptyset$.
\item[$n=2$] As seen in subsection \ref{ss_gpo}, the extreme rays for
$\cone{2}$ are 
$h^2_\emptyset = f^2_\emptyset$ and $h^2_1=f^2_1-f^2_\emptyset$.
\item[$n=3$] By repeated use of Theorems \ref{T_per} and \ref{C_select},
we may generate five extremes for $\cone{3}$: $h^3_\emptyset$, $h^3_1$,
$h^3_2$, $h^1_\emptyset \ast h^2_1 = f^3_{12} - f^3_1$ and
$h^2_1\ast h^1_\emptyset = f^3_{12} - f^3_2$.  Direct computation
shows these are all the extremes in this case.
\item[$n=4$]  Direct calculation shows that there are 13 extremes
for $\cone{4}$.  All except one can be obtained by repeated use of
Theorems \ref{T_per} and \ref{C_select}.  The remaining one is
$f^4_{13}-f^4_1+f^4_2-f^4_3$, which represents the inequality
given in Example \ref{banker}.
\item[$n=5$]  Again, direct calculation reveals that $\cone{5}$ has
41 extremes.  All but seven of these arise from lifting and convolution
as above.  The remainder are
\begin{enumerate}
\item $f^5_{134}-f^5_{14}+f^5_{24}-f^5_{34}-f^5_2+f^5_3$ 
\item $f^5_{124}-f^5_{12}+f^5_{13}-f^5_{14}+f^5_2-f^5_3$. 
\item $f^5_{1234}-f^5_{123}-f^5_{234} +f^5_{13}-f^5_{14}+f^5_{23}+f^5_{24}-f^5_2$. 
\item $f^5_{1234}-f^5_{123}-f^5_{234} 
+f^5_{13} -f^5_{14}+f^5_{23}+f^5_{24}-f^5_3$.  
\item $f^5_{124}+f^5_{234}
-f^5_{12}+f^5_{13}-f^5_{14}-f^5_{23}-f^5_{24}+f^5_2$. 
\item $f^5_{123}+f^5_{134}
-f^5_{34}+f^5_{24}-f^5_{14}-f^5_{23}-f^5_{13}+f^5_3$ .
\item $f^5_{134}+f^5_{124}-f^5_{13}-f^5_{14}+f^5_{23}-f^5_{24}$. 
\end{enumerate}
\end{description}

For rank $6$ direct calculation yields 796 extreme rays. Only 131 of 
them come from earlier extremes via lifting and convolution; 
the remaining 665 are new. An interesting problem would be to a find a reasonable characterization of the extreme rays of $\cone{n}$.


\end{document}